\documentclass[12pt,reqno,a4text]{amsart}
\usepackage{latexsym}
\usepackage{amsmath}
\usepackage{amssymb}
\usepackage[dvips]{graphicx}
\usepackage{latexsym}
\usepackage{enumerate,graphicx}
\usepackage[toc,page]{appendix}
\usepackage{setspace}
\usepackage{color,soul}
\usepackage{multirow}
\usepackage{a4wide}
\usepackage{amsthm}  
\usepackage{amsaddr}

\usepackage[colorlinks=true]{hyperref}

\newcommand{\R}{\mathbb{R}}

\def\be{\begin{equation}}
\def\ee{\end{equation}}

\def\t{\tilde}

\def\wh{\widehat}

\def\g{\gamma}

\def\p{\partial}

\def\dsl{\displaystyle}

\def\ep{\epsilon}

\def\w{\wedge}

\def\a{\alpha}

\newtheorem{theo}{Theorem}[section] 
\newtheorem{prop}[theo]{Proposition}  
\newtheorem{lemma}[theo]{Lemma}

\newtheorem{rem}[theo]{Remark}
\theoremstyle{definition}
\newtheorem{exam}{Example}[section]

\newcommand{\fp}{\frak p}

\newcommand{\fc}{\frak c}

\newcommand{\Rho}{{\mbox{\sf P}}}

\begin{document}
\title[Conformal Mercator equation \& conserved quantities of conformal circles]{Conformally flat Mercator equation and conserved quantities of conformal circles}

\author{Prim Plansangkate$^1$ and Lenka Zalabov\' a$^2$}

\address{$^1$Division of Computational Science, Faculty of Science, Prince of Songkla University, Hat Yai, Songkhla, 90110 Thailand}
\email{prim.p@psu.ac.th}

\address{$^2$Department of Mathematics,
Faculty of Science, University of South Bohemia, \\ Brani\v sovsk\' a 1760, 370 05 \v Cesk\' e Bud\v ejovice, Czechia
}
\email{lzalabova@gmail.com}

\date{}



\vskip 50pt

\begin{abstract} 
We study properties of solutions of the conformal Mercator equation on the conformal sphere. We focus on conformal circles and give two different approaches of finding their conserved quantities from conformal Killing vector fields - variational and tractorial. We find a direct relation between the two sets of quantities.
\end{abstract}

\maketitle



\section{Introduction}
Conformal circles are the most degenerate examples of conformal curves on the  conformal sphere. There are various approaches to study these curves. Tractor calculus \cite{BEG94, SZ19} is a well-established theory to study conformal curves by means of conformal curvatures and an invariant tractor Frenet-like frame. Then conformal circles are those curves for which all conformal curvatures degenerate. Other examples of such curves are loxodromes \cite{E23}, which can be viewed as having the first conformal curvature constant and all other curvatures vanishing \cite{Sulanke, Sulanke2}. We introduce these invariants and tractors in Sections \ref{st} and \ref{subsec:conservedcircle}.

The most traditional approach to studying conformal curves is to view them as solutions to suitable conformally invariant ODEs \cite{BE90,E23,DK21}. For example, conformal circles are defined as solutions of a third-order ODE \cite{Y57, BE90, T12}.  In this respect, it is natural to ask whether the equation for a class of conformal curves is variational, and to focus on conserved quantities of the solutions.  
For conformal circles, it was proved in \cite{K25} that the conformal circle equation in any dimension greater than $3$ is not variational in the classical sense.  
However, if one goes beyond the classical variation, there exists in any dimension a variational formulation for the conformal circle equation \cite[Section 3]{BE90}, associated with a conformally invariant $3$rd order Lagrangian \cite{DK21}. 

This conformally invariant Lagrangian leads to examples of variational conserved quantities of conformal circles in the conformally flat case.
It is due to the work of \cite{PZ}, where the $4$th order Euler-Lagrange equation under classical variation of the Lagrangian  was considered.  The equation, called the conformal Mercator equation, admits conformal circles as solutions in the conformally flat case \cite[Section 5]{DK21}. Therefore the conserved quantities along solutions of the conformal Mercator equation, found using Noether's theorem, in \cite[Section 3]{PZ},
are examples of conserved quantities along conformal circles that come from variational method. The quantities, which are constructed from generators of conformal symmetries of the Lagrangian, i.e. conformal Killing vector fields, are reviewed in Section 2.2.  In addition, we study properties of solutions of the conformal Mercator equation and present a characterization of solutions in terms of tractor quantities in Section 2.5.

There is another approach to construct conserved quantities of conformal circles: one can use BGG theory \cite{CG17} and tractors parallel along the curves for the tractor derivative. Such a method was introduced in \cite{GST21} for conformal circles and modified to the case of special loxodromes in \cite[Section 2.3.]{PZ}. We focus here on parallel tractors of conformal circles corresponding to conformal Killing vector fields. We find such a parallel tractor in Section \ref{subsec:conservedtractor} and give the corresponding conserved quantities in Section \ref{subsec:conservedquant}.
This allows us to study the relations between the tractor quantities and the variational quantities coming from Noether's theorem which canonically arise from conformal Killing vector fields.  Our main result, given in  Section \ref{subsec:compare}, is the identification of the two sets of conserved quantities, which shows that the conserved quantities from tractors are variational.


\section{ Conformal Mercator equation on conformal sphere} \label{sphereMercator} 
\subsection{Conformal sphere}

The conformal sphere is the sphere $S^n$ with the conformal class of metrics represented by the standard round metric. We can view it as the projectivization of the cone of non-zero null vectors in the pseudo-Euclidean space $\R^{n+1,1}$ of signature $(n+1,1)$.
For explicit computations in coordinates, we represent the corresponding standard scalar product 
by the block $(1,n,1)$-matrix
\begin{align} \label{met}
G=\left( \begin{smallmatrix}
0&0&1 \\ 0& {\rm Id} &0 \\ 1 &0&0
 \end{smallmatrix} \right)
\end{align}
for the identity matrix ${\rm Id}$.  
The conformal sphere is homogeneous, and the Lie group of conformal transformations of $S^n$ is isomorphic to $O(n+1,1),$ with the corresponding Lie algebra $\frak{so}(n+1,1).$  
The stabilizer of a point is the Poincar\'e subgroup $P\subset O(n+1,1)$ and $S^n\cong O(n+1,1)/P$.
According to the metric \eqref{met}, elements of $\frak{so}(n+1,1)$  are the block $(1,n,1)$-matrices
\begin{align} \label{sym-matrix}
\left( \begin{smallmatrix} -a&2S&0 \\ T &-R&-2S^t \\ 0&-T^t&a \end{smallmatrix} \right)
\end{align}
for $a \in \R$, $R \in \frak{so}(n)$ and $T,S^t \in \R^n$, 
and the Poincar\' e conformal algebra $\frak{p} \subset \frak{so}(n+1,1)$ corresponds to the upper triangular matrices.

We get local exponential coordinates on $S^n$ around the origin $o=eP$ corresponding to the usual coordinates coming from the projectivization of the cone in $\R^{n+1,1}$ as follows:
We fix the complement $\fc$ of $\frak{p}$ in $\frak{so}(n+1,1)$ such that $\frak{so}(n+1,1)=\fc \oplus \fp$ of the form
\begin{align}   \label{decomp}
\fp \ni \left( \begin{smallmatrix} -a&2S&0 \\ 0 &-R&-2S^t \\ 0&0&a \end{smallmatrix} \right), \ \ 
\fc \ni \left( \begin{smallmatrix} 0&0&0 \\ T &0&0 \\ 0&-T^t&0 \end{smallmatrix} \right).
\end{align}
Thus $\fc \simeq \R^n$ and denoting $\langle\: , \: \rangle$  the standard scalar product on $\R^n$, we compute
\begin{align*}
\exp \left( \begin{smallmatrix} 0&0&0 \\ T &0&0 \\ 0&-T^t&0 \end{smallmatrix} \right) \left( \begin{smallmatrix} 1 \\ 0 \\ 0 \end{smallmatrix} \right) = \left( \begin{smallmatrix} 1&0&0 \\ T & {\rm Id} &0 \\ -\frac12\langle T,T \rangle &-T^t&1 \end{smallmatrix} \right)  \left( \begin{smallmatrix} 1 \\ 0 \\ 0 \end{smallmatrix} \right) = \left( \begin{smallmatrix} 1 \\ T\\ -\frac12 \langle T,T \rangle \end{smallmatrix} \right).
\end{align*}
Then, from now on we can treat everything on the sphere as a local computation on $\R^n$.
Thus the elements of \eqref{sym-matrix}, viewed as conformal Killing vector fields, take the coordinate form
\begin{align} \label{kf}
(T^i-{R^i}_j x^j+ax^i-2S_jx^jx^i+x_jx^jS^i)\frac{\partial}{\partial x^i}.
\end{align}

\subsection{Conformal Mercator equation}
In local coordinates, let $\gamma=(\gamma^i(t))$, $i=1, \dots, n$, be an arbitrary curve 
 and 
$$U=(\frac{d}{dt}\gamma(t)), \ \ A=(\frac{d^2}{dt^2}\gamma(t)), \ \ A'=(\frac{d^3}{dt^3}\gamma(t)), \dots$$
be its velocity vector field, acceleration vector field, jerk vector field and so on for subsequent derivatives, where $\,{}'$ denotes $\frac{d}{dt}.$
We denote by $u$ the length of the velocity vector field
$$u=\sqrt{\langle U,U \rangle}.$$
The conformal Mercator equation \cite{DK21} on the conformal sphere is given by 
 \begin{align}
 \begin{split}
     \label{confMeqflat} \frac{dC}{dt} &= 0, \\
 C &= u^{-2} \left( A' - u^{-2}|A|^2 U -2 u^{-2} \langle A,U\rangle A + 4 u^{-4} \langle A,U \rangle^2 U - 2 u^{-2} \langle A',U \rangle U \right).
  \end{split}
 \end{align}  
It is the Euler-Lagrange equation of the Lagrangian \cite[Section 4.2]{DK21},
\[L =  \frac{d}{dt} \left( u^{-2}\langle U,A \rangle \right) + \frac{1}{2} u^{-2}|A|^2 - u^{-4} \langle U,A \rangle^2,\]
which is conformally invariant.  That is, for any one-parameter group of transformations $\{\phi_s\}$ generated by a conformal Killing vector field, 
the Lagrangian $L$ is preserved in the sense that $L(\g(t)) = L( \phi_s(\g(t))),$
where $L(\g(t))$ denotes the Lagrangian evaluated on a curve $\g,$ and $L( \phi_s(\g(t)))$ the Lagrangian evaluated on the transformed curve $\phi_s(\g(t)).$
 Then, Noether's theorem gives a conserved quantity (first integral) along solutions of the conformal Mercator equation \eqref{confMeqflat}, for each symmetry group.  Such conserved quantities are constructed in \cite{PZ}.
In particular, it was proved \cite[Proposition 3.1]{PZ} that for any solution $\g$ of \eqref{confMeqflat} and any conformal Killing vector field $V$ restricted to $\g,$  the function
\be \label{F} F = \frac{d}{dt} \langle W, V' \rangle + \langle W', V' \rangle  - \langle C,V \rangle,  \qquad W = u^{-2} U,\ee
where $C$ is given in \eqref{confMeqflat}, is constant along $\g.$

To discuss properties of conserved quantities of conformal circles later in Section \ref{par-trac}, we shall make use of the expressions of the conserved quantities of the form \eqref{F}, where $V$   corresponds to a vector field in a specific basis of the Lie algebra of conformal Killing vector fields.  These are  given in the following proposition.
Below we use the elementary alternating tensor field to express the quantities.  Recall that the elementary alternating tensor field  of rank $m,$ denoted by $\ep^{i_1 \dots i_m},$ is defined by its action on $m$ vector fields, $Y, \dots, Z,$ as
\[ \ep^{i_1 \dots i_m} (Y, \dots, Z) 
= 
\det \left(
\begin{array}{ccc}
Y^{i_1} & \dots &  Z^{i_1}  \\
\vdots & { } &  \vdots \\
Y^{i_m}& \dots &  Z^{i_m}
\end{array}
\right),
\]
where $Y^j$ and $Z^j$ are components of $Y$ and $Z,$ respectively.

\begin{prop} \cite[Proposition 3.2]{PZ} \label{prop_FTFRFDFS}
Let $(x^i)$ be the standard coordinates for $\R^n,$ and $\g$ be a solution to the conformal Mercator equation {\rm \eqref{confMeqflat}} given by $(x^i) = (\g^i(t)) =: (X^i(t)),$ with $U$ and $A$ denoting its velocity and acceleration vector fields, respectively. Let $X$ denote the vector field $X^i  \frac{\p}{\p x^i}$ along $\g.$ Then, the following quantities are constant on $\g:$
\begin{eqnarray}
E_T^i  &=& - C^i \label{FT} \\
E_R^{ij} &=& u^{-2}\ep^{ij}\left(U, A\right) + \ep^{ij}\left(C, X\right)\label{FR} \\
E_D &=& -u^{-2} \langle U,A \rangle - \langle C,X \rangle \label{FD} \\
E_S^i &=& 2 u^{-2} \langle U, X \rangle A^i - 2 \left(1+ u^{-2} \langle A, X \rangle \right) U^i \nonumber \\
&{}& + 2 \left(u^{-2} \langle U,A \rangle + \langle C,X \rangle \right) X^i - |X|^2 C^i,  \label{FS}
\end{eqnarray}
where $C$ is given by \eqref{confMeqflat}. 

The functions \eqref{FT} - \eqref{FS} are obtained from \eqref{F}, where $V$ is given by \eqref{kf} with specific values of $T^i,$ $R^i_j,$ $a$ and $S^i.$
In particular, $E_T^i,$ $E_R^{ij},$ $E_D$ and $E_S^i$ are derived from the generators of translation $\frac{\p}{\p x^i}$, rotation $x^i \frac{\p}{\p x^j} - x^j \frac{\p}{\p x^i}$, dilatation $x^i \frac{\p}{\p x^i}$ and special conformal transformation $x_j x^j  \frac{\p}{\p x^i} - 2 x_i x^j \frac{\p}{\p x^j},$ respectively. 
\end{prop}

\subsection{Standard tractors and invariants of curves} \label{st}

To discuss properties of solutions of the Mercator equation in detail, we also employ a tractor viewpoint on curves; see \cite{CG17} for an introduction to tractor calculus. Due to the homogeneity of the conformal sphere, we follow \cite{GZ} to describe tractors in coordinates. 

We consider the standard tractor bundle $\mathcal{T}$ modelled on $\R^{n+2}$ with the tractor metric given by \eqref{met}. Thus, we deal with the standard representation $$\rho: \frak{so}(n+1,1) \to \frak{gl}(\R^{n+2}).$$ 
The tractors are then given in the above exponential coordinates just by restriction to $\fc$.  The standard tractor bundle $\mathcal{T}$ trivialises as $\mathcal{T} \simeq \fc \times \R^{n+2}$ and  standard tractors $\mathbb{X} \in \mathcal{T}$ are sections of $\fc \times \R^{n+2}$, that is, functions $ \fc \to \R^{n+2}$. 
In this local trivialization, the corresponding functions $\mathbb{X}$ take the form
\begin{gather} \label{exp}
\exp(J) \mapsto \exp(-\rho(J))(w), \ \ \ \ J=\left( \begin{smallmatrix} 0&0&0 \\ X &0&0 \\ 0&-X^t&0 \end{smallmatrix} \right)
\end{gather}
for some neighborhood of $0$ in $\fc$, i.e. $X \in \fc$ and $w=(w^a)$ are coordinates in the representation space $\R^{n+2}$.

For each curve $\gamma$, we can define a canonical series of derived standard tractors by means of the tractor derivative. The standard tractor connection $\nabla^\mathcal{T}$ decomposes into the fundamental derivative $D$ and the algebraic action, and in the case of the conformal sphere takes the form 
$\nabla^\mathcal{T}=D+\rho$ restricted to $\fc$ to get coordinate description.
Unlike the usual way of writing standard tractors in columns, \cite{CG17,SZ19,GST21}, we write them in coordinates in the standard basis $(e_a)$ of $\R^{n+2}$, where $a=0, \dots, N$ and $N:=n+1$. For later purposes, this notation is more useful for dealing with higher powers of the representation $\rho$. 

We decompose the standard tractor bundle by means of the representation space for the canonical flat metric from the conformal class as 
$$ \R^{1+n+1}=\langle e_0 \rangle + \langle e_1, \dots, e_n \rangle + \langle e_{N} \rangle.
$$
Here $e_0$ corresponds to the injective slot and $e_N$ to the projective slot of the standard tractor.
We start the series of standard tractors with the canonical inclusion 
\begin{gather} \label{trT}
\mathbb{T}= u^{-1}e_0.
\end{gather}
Then we compute the tractor derivative $U^a\nabla^\mathcal{T}_a$ along the curve $\gamma$, which takes the explicit form
\begin{align} \label{tra}
U^a\nabla^\mathcal{T}_a(\mathbb{X})=D(\mathbb{X})+ \rho(K)(\mathbb{X}),\ \ \ \  K:=\left( \begin{smallmatrix} 0&0&0 \\ U&0&0 \\ 0&-U&0 \end{smallmatrix} \right)
\end{align}
for the standard tractor $\mathbb{X}$.
 Direct computation using \eqref{tra} leads to the sequence of standard tractors $\mathbb{U}= U^a\nabla^\mathcal{T}_a(\mathbb{T})$, $\mathbb{A}= U^a\nabla^\mathcal{T}_a(\mathbb{U})$, $\mathbb{A'}= U^a\nabla^\mathcal{T}_a(\mathbb{A})$ etc. as follows,
 where  $i=1,\dots,n$  (compare e.g. with \cite[Section 2.1.]{SZ19} for ${\sc \Rho}=0$)
\begin{align}
\mathbb{U}&=-u^{-3}\langle U,A \rangle e_0+u^{-1} U^ie_i \label{trU}
\\
\mathbb{A}&= (3u^{-5}\langle U,A \rangle^2-u^{-3}(\langle A,A \rangle +\langle U,A'\rangle ))e_0+
(-2u^{-3} \langle U,A \rangle U^i+u^{-1}A^i)e_i-u e_{N} \label{trA}
\\
\begin{split}
\mathbb{A}'&= (-15u^{-7}\langle U,A \rangle^3+9u^{-5}\langle U,A \rangle(\langle A,A \rangle +\langle U,A'\rangle)
-u^{-3}(3\langle A,A' \rangle +\langle U,A'' \rangle) )e_0
\\ &+
(9u^{-5}\langle U,A \rangle^2U^i-3u^{-3}(\langle A,A \rangle +\langle U,A'\rangle)U^i-3u^{-3}\langle U,A \rangle A^i+u^{-1}(A')^i)e_i
\label{trAbar}
\end{split}
 \\
 &\textrm{etc.} \notag
\end{align}

For each (parametrized) curve, thus for each sequence of tractors $\mathbb{T}, \mathbb{U},  \mathbb{A}, \mathbb{A'}, \dots $ of arbitrary length $\ell$,
we compute the following Gram matrix of order $\ell$ 
 \begin{align} \label{det}
M_\ell=\left( \begin{smallmatrix}
  0&0&-1&0&\alpha_1 & \cdots\\
   0&1&0&-\alpha_1&-\frac32 \alpha'_1
 & \cdots\\
   -1&0&\alpha_1&\frac12\alpha'_1&\frac12 \alpha''_1-\alpha_2& \cdots\\
  0&-\alpha_1&\frac12 \alpha'_1&\alpha_2&\frac12 \alpha'_2& \cdots\\
   \alpha_1&-\frac32 \alpha'_1&\frac12 \alpha''_1-\alpha_2&\frac12 \alpha'_2&\alpha_3 & \cdots \\
   \vdots &  \vdots & \vdots & \vdots & \vdots & \ddots
  \end{smallmatrix} \right)
\end{align}
for the standard tractor metric, compare to \cite[formula (13)]{SZ19}. Thus $\alpha_1$ is the square of the length of the acceleration tractor  $\mathbb{A}$ and analogously for next tractors/alphas in the series.
In particular, we compute
\begin{align} \label{alpha1}
\alpha_1 &= -6u^{-4}\langle U,A \rangle^2+2u^{-2}\langle U, A' \rangle +3u^{-2}\langle A,A \rangle, \\ 
\begin{split}
\label{alpha2}
\alpha_2 &=  27 u^{-8}\langle U,A \rangle^4-18 u^{-6}\langle U,A \rangle^2\langle U,A' \rangle -27 u^{-6}\langle U,A \rangle^2\langle A,A \rangle +3u^{-4}\langle U,A' \rangle^2 \\ &-6u^{-4} \langle U,A \rangle\langle A,A' \rangle+12u^{-4} \langle U,A' \rangle \langle A,A \rangle +9u^{-4} \langle A,A \rangle^2 +u^{-2} \langle A',A' \rangle
\end{split}
\end{align}

This also allows us to define a series of natural relative invariants $\Delta_\ell=\det(M_\ell)$, \cite[Lemma 2.2]{SZ19}.
Thus, the invariants $\Delta_\ell$ depend on the parametrisation; however, their vanishing, respectively non-vanishing, is independent of reparametrizations.
Let us note that $\Delta_1=\Delta_2=0$ and $\Delta_3=-1$. In particular, the tractors $\mathbb{T}$, $\mathbb{U}$ and $\mathbb{A}$ are linearly independent.
The first non-trivial invariant then is $\Delta_4= \alpha_1^2-\alpha_2$, and  
\begin{align}
\begin{split} \label{Delta}
\Delta_4&=
9u^{-8}\langle U,A \rangle^4-6u^{-6}\langle U,A \rangle^2\langle U,A'\rangle-9u^{-6}\langle U,A \rangle^2\langle A,A\rangle \\ 
&+6u^{-4}\langle U,A\rangle\langle A,A'\rangle +u^{-4}\langle U, A' \rangle^2-u^{-2}\langle A',A'\rangle.
\end{split}
\end{align}

 Finally, let us note that there is a series of absolute conformal invariants, aka conformal curvatures, that take the form 
 \begin{align}
\kappa_1 &= -\frac12 (-\Delta_4)^{-\frac52} (\alpha_1 \Delta_4^2-\frac12 \Delta_4\Delta_4''+\frac{1}{16}(\Delta_4')^2), \label{kappa1}\\
\kappa_i &= -(-\Delta_4)^{-\frac14}(\Delta_{i+1} \Delta_{i+3})^\frac12 \Delta_{i+2}^{-1},\ \  i>1
 \end{align} 
 These are independent of the choice of parametrisation.
The first curvature is exceptional. In particular, it is not defined for $\Delta_4=0$, which is a condition that distinguishes conformal circles. According to \cite{Sulanke,Sulanke2, SZ19}, conformal circles are the most degenerate in the sense that all conformal curvatures are undefined. The second curvature takes the form
 \begin{align}
\kappa_2 &= (-\Delta_4)^{-\frac54}(- \Delta_{5})^\frac12 
 \end{align}
Thus $\kappa_2$ vanishes for $\Delta_5=0$ if defined, i.e. if $\Delta_4 \neq 0$. According to \cite{Sulanke,Sulanke2, E23}, if such a curve moreover satisfies $\kappa_1$ constant, we call it a loxodrome. This is particularly satisfied for 
 curves such that $\Delta_5=0, \Delta_4 \neq 0, \alpha_1'=0, \Delta_4'=0$.

\begin{rem}
R. Sulanke in his work \cite{Sulanke, Sulanke2} also finds further curves with constant conformal curvature.
On the conformal sphere $S^n$ for $n\geq 3$, generic curves are determined by so-called conformal arc-length (related to $\Delta_4$) and $n-1$ conformal curvatures $\kappa_i$. Conformal circles and loxodromes are the simplest ones.
Conformal circles are precisely curves contained in a fixed round circle $S^1 \subset S^{n}$, and loxodromes lie in a conformal $2$-sphere $S^2\subset S^{n}$.
\begin{center}
\begin{figure}[h!]
    \includegraphics[width=0.3\linewidth]{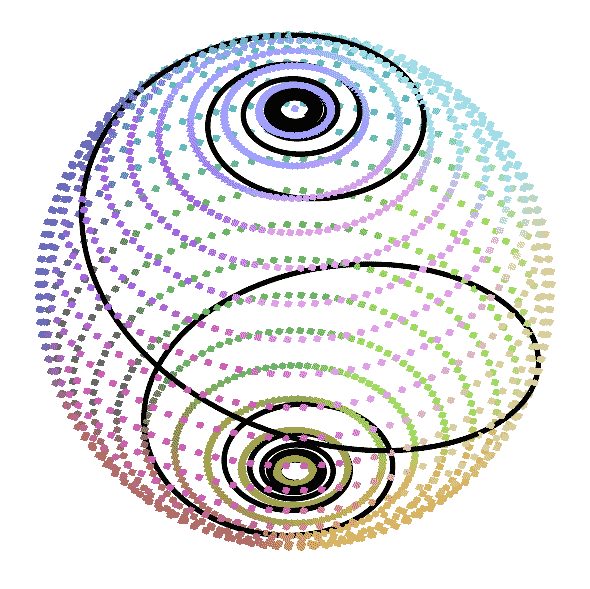}
    \caption{Family of conformal circles (dotted line) and a loxodrome (solid line) on $S^2$.}
    \label{curves}
\end{figure}
\end{center}
\end{rem}


\subsection{Tractor equations} \label{subsec:conservedcircle}
There is a natural question about properties of curves given by the linear dependency of the series of the tractors $\mathbb{T}, \mathbb{U}, \mathbb{A}, \mathbb{A'}, \dots, \mathbb{A}^{\ell}$, i.e. vanishing of $\Delta_{\ell+3}$.
The first non-trivial case is the case of $\ell=1$.

Recall that a curve is an (arbitrarily parametrised) conformal circle if and only if 
the  tractors $\mathbb{T}, \mathbb{U}, \mathbb{A}, \mathbb{A'}$
are linearly dependent. Then
$\mathbb{A'} \in \langle \mathbb{T}, \mathbb{U}, \mathbb{A} \rangle$ and it follows that the curve is a conformal circle if and only if $\Delta_4=0$, \cite[Proposition 3.3]{SZ19}. One finds
\begin{align} \mathbb{A'}=-\alpha_1\mathbb{U}-\frac12 \alpha'_1 \mathbb{T}  
\label{tc}\end{align}
using the Cramer rule and \eqref{det}.
In the following, the condition $\alpha_1=k$ for a constant $k$ will be important.
Involving this condition, we get the tractor equation 
\be \label{tractor-eq-k} \mathbb{A'}+k\mathbb{U}=0. \ee 

The middle (conformal invariant) slot of \eqref{tc} 
then gives an ODE for conformal circles, which takes the general form
$$
u^{-1}A'-3u^{-3}\langle U,A\rangle A+3u^{-5}\langle U,A\rangle^2 U-u^{-3}\langle U,A'\rangle U=0.
$$
The condition $\a_1 =k$ can be written as 
\begin{align} \label{subk}
\dsl \langle U, A' \rangle = 3u^{-2}\langle U,A \rangle^2-\frac32 \langle A,A \rangle +\frac{u^2k}{2}
\end{align}
according to \eqref{alpha1}. 
Thus, in the case $\alpha_1=k$ we get the  equation, i.e. the middle (invariant) slot of \eqref{tractor-eq-k} as
\begin{align}
\label{equationk}
u^{-1}A'-3u^{-3} \langle U,A \rangle A+\frac32u^{-3} \langle A,A \rangle U-\frac12 u^{-1}kU=0.
\end{align}

\begin{rem} 
\;
\begin{enumerate}[{\rm(1)}]
\item Note that the tractor $\mathbb{T}$ has a zero middle slot and does not contribute to the ODEs.
\item In the case $k=0$, the equation \eqref{equationk} becomes the usual equation of (projectively parametrised) conformal circles
(up to the multiplication of a suitable power of $u$), \cite{BE90}.
\end{enumerate}
\end{rem}


Let us now consider the case $\ell=2$, i.e. vanishing of $\Delta_5$. Thus we assume $\Delta_4 \neq 0$ and $\Delta_5=0$, which means that the tractors $\mathbb{T}, \mathbb{U},  \mathbb{A},  \mathbb{A}', \mathbb{A}''$ are linearly dependent, i.e. $ \mathbb{A}'' \in \langle \mathbb{T}, \mathbb{U},  \mathbb{A},  \mathbb{A}'\rangle$, \cite{PZ}. We again find the tractor combination explicitly using \eqref{det} and the Cramer rule as
\begin{align} \begin{split}
\label{tract-comb-5}
\mathbb{A}''&=\frac{\Delta'_4}{2\Delta_4}\mathbb{A}'-\alpha_1\mathbb{A}+\frac{1}{2\Delta_4}(\alpha'_1(2\alpha_2-\Delta_4)-\alpha_1\alpha'_2)\mathbb{U}\\ &+
\frac{1}{4\Delta_4}(2\alpha_1(\alpha'_1)^2-4\Delta_4^2-2\alpha''_2\Delta_4-\alpha'_1\alpha'_2)\mathbb{T}.
\end{split}
\end{align}
We do not write the invariant slot, i.e., the ODEs, for the corresponding curves in full generality, since we do not need it and it is very long.  For what follows, we are only interested in the case $\Delta_4'=0$ (and $\Delta_4 \ne0$). With this assumption, the coefficient of $\mathbb{A'}$ vanishes and we get the tractor equation
 \begin{align} \begin{split}
 \label{tract-comb-5-Dconst}
\mathbb{A}''&=-\alpha_1\mathbb{A}+\frac{1}{2\Delta_4}(\alpha'_1(2\alpha_2-\Delta_4)-\alpha_1\alpha'_2)\mathbb{U}\\ &+
\frac{1}{4\Delta_4}(2\alpha_1(\alpha'_1)^2-4\Delta_4^2-2\alpha''_2\Delta_4-\alpha'_1\alpha'_2)\mathbb{T}.
\end{split}
\end{align}
The middle invariant slot, i.e. ODEs for curves such that $\Delta_5=0$, $\Delta_4 \neq 0$, $\Delta_4'=0$ take the form (compare \cite[Appendix C]{PZ})
\begin{align} 
 \begin{split}
 \label{tr}
E:= &(12 u^{-7}\langle U, A \rangle^3 -8 u^{-5} \langle U,A \rangle \langle U, A' \rangle-3u^{-5}\langle U,A \rangle \langle A,A \rangle +u^{-3}\langle U,A'' \rangle)U+\\
&(-12u^{-5}\langle U,A\rangle^{2}+4u^{-3}\langle U, A' \rangle+3u^{-3}\langle A, A \rangle)A +4u^{-3}\langle U,A\rangle A' -u^{-1}A''=0.
\end{split}
 \end{align}


\subsection{Two types of solutions of the conformal Mercator equation} \label{sec:MercSol}
 
Let us return to the conformal Mercator equation.
First, we shall observe some properties of solutions of the equation in the following lemma.

\begin{lemma} \label{lem_DalMer}
Any solution of the conformal Mercator equation \eqref{confMeqflat} satisfies 
$$\Delta_4'=0\ \ \  \rm{and} \ \ \ \a_1'=0.$$
\end{lemma}

\begin{proof}
First, one can write $\Delta_4$ \eqref{Delta} in terms of the vector $C$ defined in the conformal Mercator equation \eqref{confMeqflat} as
\[ \Delta_4 = -u^2|C|^2+ \big(\langle C, U \rangle - u^{-4} \langle U, A \rangle^2 \big)^2 + 2 u^{-2} \langle U, A \rangle \langle C, A \rangle - u^{-6} \langle U, A \rangle^2 |A|^2.\]
Hence, 
\be \label{dD4C} \frac{1}{2}\Delta_4' = u^{-2} \langle \frac{dC}{dt}, A \rangle \langle U, A \rangle +   \langle \frac{dC}{dt}, U \rangle  \big(\langle C, U \rangle - u^{-4} \langle U, A \rangle^2 \big) - u^2 \langle \frac{dC}{dt}, C \rangle. \ee

Also, a direct calculation shows a relation between the derivative $\frac{dC}{dt}$  and the expression $E$ on the left-hand side of \eqref{tr}  given by
\be \label{connectMwd5} \frac{dC}{dt} = - u^{-1} E - \frac{1}{2}u^{-2} \a_1' U.\ee
In addition, it can be shown that $\langle E, U\rangle \equiv 0,$ which 
implies that 
\be \label{CpUal1} \langle \frac{dC}{dt}, U \rangle = - \frac{1}{2} \a_1'. \ee
Equations \eqref{dD4C} and \eqref{CpUal1} then show that any solution of the conformal Mercator equation, $\frac{dC}{dt}=0,$ satisfies $\Delta_4'=0$ and $\a_1'=0.$ 
\end{proof}

Next, the following proposition characterises all solutions of the conformal Mercator equation in terms of the relative invariants $\Delta_\ell$ and the condition on $\a_1.$

\begin{prop}
Solutions of the conformal Mercator equation \eqref{confMeqflat} consist of two types of curves; 

Type I: Curves satisfying $\Delta_4=0$ and $\a_1'=0,$ and

Type II: Curves satisfying $\Delta_5 =0,$ $\Delta_4 \ne 0,$ $\Delta_4'=0$ and $\a_1'=0.$ 
\end{prop}

\begin{proof}
 Any curves with $\Delta_4=0$ and $\a_1'=0$ obeys equation \eqref{equationk}, the conformal circle equation with constant $\a_1 = k.$  Assuming this, the vector $C$  in \eqref{confMeqflat} reduces to 
\be \label{Ck} C= \left(\frac{1}{2} u^{-4}|A|^2 - 2 u^{-6} \langle A,U \rangle^2 - \frac{1}{2}u^{-2}k\right)U + u^{-4} \langle A,U \rangle A.\ee
Then by direct computation, it can be shown that $\frac{dC}{dt} = 0$ along any solution of \eqref{equationk}, which means any solution of \eqref{equationk} satisfies  \eqref{confMeqflat}.   We call such solutions type  I.

\vspace{0.1cm}

Next, by Lemma \ref{lem_DalMer} the complementary of solutions of type I are those with nonzero constant $\Delta_4$ and $\a_1'=0.$
The equation for the curves with $\Delta_5 =0,$ $\Delta_4 \ne 0,$ $\Delta_4' = 0$ is given by \eqref{tr}.  Then from \eqref{connectMwd5},
\[\frac{dC}{dt} = - u^{-1} E - \frac{1}{2}u^{-2} \a_1' U,\]
it follows that any solution of the conformal Mercator equation, satisfying $\frac{dC}{dt} = 0$ and $\a_1'=0,$ obeys $E=0,$ i.e. \eqref{tr}, and any solution of \eqref{tr} with $\a_1'=0$ also solves the conformal Mercator equation \eqref{confMeqflat}.  This shows that the remaining solutions, named type II, are all curves satisfying $\Delta_5 =0,$ $\Delta_4 \ne 0,$ $\Delta_4'=0$ and $\a_1'=0.$  
\end{proof}

\begin{rem} \;
    \begin{enumerate}[{\rm(1)}]
        \item It was already noted in \cite{DK21} that projectively parametrized conformal circles (curves with $\Delta_4=0, \a_1=0$) satisfy the conformal Mercator equation \eqref{confMeqflat}.  Our proof for solutions of type I uses the same argument adapted to the condition $\a_1'=0.$  
        \item Although \eqref{tr} is derived under the assumption of $\Delta_4 \ne 0,$ it can be directly verified that any solution of \eqref{equationk}, i.e. a parametrized conformal circle with constant $\a_1=k,$ satisfies \eqref{tr}.
    \end{enumerate}
\end{rem}

\vspace{0.2cm}

The next proposition says what kind of traces or trajectories (unparametrized curves) the solutions of the conformal Mercator equation \eqref{confMeqflat} represent.

\begin{prop}
Consider solutions of the conformal Mercator equation \eqref{confMeqflat}.
The traces of solutions of type I cover all conformal circles, and the traces of type II cover all loxodromes.
\end{prop}

\begin{proof}
With $\Delta_4=0,$ it is clear that the traces of solutions of type I are conformal circles, and any conformal circle (in fact any curve) admits a reparametrization so that $\a_1$ is constant. To see this, recall that under a reparametrization $\t t = g(t),$ $\a_1$ transforms to
$$\t \a_1 = g'^{-2}(\a_1-2S(g)),
\ \ {\rm where}\ \  S(g) =\frac{g'''}{g'} -\frac{3}{2} \left(\frac{g''}{g'}\right)^2.$$ Then $\t \a_1 = k$ constant is an ODE, which always admits a solution at least locally. Cf. \cite[Section 4]{BE90}.

 For type II, recall that loxodromes are curves satisfying $\Delta_5 = 0,$ $\Delta_4 \ne0$ and $\kappa_1$ constant. Hence the traces of solutions of type II are loxodromes.  Now we shall argue that any loxodrome admits a reparametrization such that it satisfies  $\Delta_5 = 0,$ $\Delta_4 \ne0,$ $\Delta_4'=0$ and $\a_1'=0.$ 

Under a reparametrization $\t t = g(t),$ $\Delta_4$ transforms to $\t \Delta_4 = g'^{-4} \Delta_4,$
and thus
given any $\Delta_4(t) \ne 0,$ 
one can always find $g(t)$ so that  $\t \Delta_4$ is constant. Therefore any curve satisfying $\Delta_5 =0,$ $\Delta_4 \ne 0$  
can be reparametrized so that it has constant ${\Delta_4}.$  That is, the condition  $\Delta_4'=0$ can be regarded as fixing 
 a preferred family of parametrizations, with the remaining freedom given by the affine transformation $g(t) = at+b,$ $a \ne 0.$  Cf. $${\t \Delta_4}' = -4 g'^{-6} g'' \Delta_4 + g'^{-5} \Delta_4'.$$

 Now, since $\kappa_1$ is an absolute invariant, i.e. it is invariant under reparametrization, any loxodrome in a preferred parametrization with $\Delta_4'=0$ has constant $\a_1,$ c.f. \eqref{kappa1}. 
 \end{proof}

 \begin{rem}
It was stated in \cite{DK21} that solutions of the conformal Mercator equation that are not conformal circles (i.e. our solutions of type II) can be obtained from conformal transformations of the logarithmic spiral on a plane,
 \[\wh X(t) = e^t  \cos(ct) P_0 +  e^t  \sin(ct) Q_0 + R_0,\]
where $c$ is a constant, $P_0, Q_0, R_0$ are constant vectors, $\langle P_0, Q_0 \rangle =0$ and ${|P_0| = |Q_0|}.$  Note also that the inverse image of such a logarithmic spiral is a loxodrome on $S^2.$
 \end{rem}

\section{Conserved quantities of conformal circles from tractors} \label{par-trac}

The BGG theory states that on the conformal sphere, each conformal Killing vector field corresponds (via BGG splitting and tractor metric) to an element of the tractor bundle for the representation 
\begin{align} \label{repr}
\rho_2:\frak{so}(n+1,1) \to \frak{gl}(\wedge^2\R^{n+2}).
\end{align}
Indeed, each such tractor $\mathbb{V}$ is parallel w.r.t. the tractor connection, because the whole tractor space corresponds to solutions of the corresponding BGG operators, and they are all normal, i.e. coming from sections parallel to the tractor connection, \cite{CGH12,GZ}.

Assume we have a tractor $\mathbb{Y}$ that is parallel along a curve $\gamma(t)$ on a conformal sphere w.r.t. the tractor connection. Then we can use these tractors $\mathbb{V},\mathbb{Y}$ to construct conserved quantities of the curve $\gamma(t)$. The tractor metric $G$ is parallel w.r.t. the tractor connection, too.  Then the pairing $G(\mathbb{V},\mathbb{Y})$ gives a conserved quantity of the curve $\gamma$. We focus here on conserved quantities constructed in this way for conformal circles.

\subsection{Parallel $2$-tractor} \label{subsec:conservedtractor}
Assume we have conformal circles satisfying the condition that $\alpha_1=k$ is constant. 
\begin{prop}
The curve $\gamma$ is a conformal circle such that $\alpha_1$ is constant if and only if the tractor $$\mathbb{Y}:=\mathbb{U} \wedge \mathbb{A}$$ is parallel for the tractor connection  $\nabla^{\wedge^2 \mathcal{T}}$ on $\wedge^2 \mathcal{T}$ along the curve $\gamma$. 
\end{prop}
\begin{proof}
Assume $\gamma$ is a conformal circle with $\alpha_1=k$ constant.
The condition $\mathbb{A'}=-k\mathbb{U}$ \eqref{tractor-eq-k} then gives 
\begin{align*}
\nabla^{\wedge^2 \mathcal{T}}(\mathbb{U} \wedge \mathbb{A}) &=
\mathbb{A} \wedge \mathbb{A} +\mathbb{U} \wedge \mathbb{A'}
\\&=\mathbb{A} \wedge \mathbb{A} +\mathbb{U} \wedge (-k\mathbb{U})
\\ &=\mathbb{A} \wedge \mathbb{A} -k (\mathbb{U} \wedge \mathbb{U})= 0.
\end{align*}
for the tractor connection $\nabla^{\wedge^2 \mathcal{T}}$ along $\gamma$. Conversely, assume $\mathbb{Y}$ is parallel for the tractor connection. Then we have
\begin{align*}
\nabla^{\wedge^2 \mathcal{T}}(\mathbb{U} \wedge \mathbb{A}) &=
\mathbb{A} \wedge \mathbb{A} +\mathbb{U} \wedge \mathbb{A'} =\mathbb{U} \wedge \mathbb{A'}
= 0.
\end{align*}
This means that for the tractors $\mathbb{U}$ and $\mathbb{A'}$, one tractor is a constant multiple of the other (or $\mathbb{A'}$ vanishes, if the multiple is zero). This means that tractors $\mathbb{T}, \mathbb{U}, \mathbb{A}, \mathbb{A'}$
are linearly dependent and we get a conformal circle.
\end{proof}
We find the explicit form of the $2$-tractor $\mathbb{Y} = \mathbb{U} \wedge \mathbb{A}$ in coordinates. 
Consider the basis 
\begin{align} \label{baze} b_2 =\{ e_{ab}=e_a \wedge e_b, \ \ \ \ 0\leq a< b \leq N\} \ \ \ {\rm of} \ \ \wedge^2\R^{n+2}.\end{align} 

 \begin{prop} \label{2tractprop}
 The $2$-tractor $\mathbb{Y} =  \mathbb{U} \wedge \mathbb{A} $ takes 
 in the case $\alpha_1=k$ the explicit form
\begin{align} \begin{split} \label{2tractor}
\mathbb{U} \wedge \mathbb{A} &= \big( (2u^{-6} \langle U,A \rangle^2   -
\frac12 u^{-4} \langle A,A \rangle   +\frac12 u^{-2} k )U^i - u^{-4} \langle U,A \rangle A^i \big) e_{0i} \\ &+u^{-2} \langle U,A \rangle e_{0N}  +u^{-2}U^{[i}A^{j]} e_{ij}-U^i e_{iN}.
\end{split} \end{align}
where we use Einstein summation convention,  $[ab..c]$ denotes the anti-symmetrization over indices in the bracket, $\langle\ , \ \rangle$ denotes the standard scalar product on $\R^n$ and $i,j=1, \dots, n$.
 \end{prop}

\begin{proof}
We proceed by direct computation with tractors \eqref{trU} and \eqref{trA}, where we involve the condition \eqref{subk} in the formulae. Then the tractors take the form
\begin{align*}
\mathbb{U}&=-u^{-3}\langle U,A \rangle e_0+u^{-1} U^ie_i 
\\
\mathbb{A}&= (
\frac12u^{-3}\langle A,A \rangle-\frac12u^{-1}k
)e_0+
(-2u^{-3} \langle U,A \rangle U^i+u^{-1}A^i)e_i-u e_{N} 
\end{align*}
We immediately omit the term with $e_0 \wedge e_0$ as it vanishes. 
We have
\begin{align*}
&-u^{-3}\langle U,A \rangle e_0 \wedge \big((-2u^{-3} \langle U,A \rangle U^i+u^{-1}A^i)e_i-u e_{N}\big ) +
\\ & u^{-1} U^ie_i \wedge \big( 
(\frac12u^{-3}\langle A,A \rangle-\frac12u^{-1}k
)e_0+(-2u^{-3} \langle U,A \rangle U^i+u^{-1}A^i)e_i-u e_{N} \big).
\end{align*}
Expanding the expression and using $U \wedge U=0$, we get 
\begin{align*}
 &2u^{-6} \langle U,A \rangle^2U^i e_0 \wedge e_i - u^{-4} \langle U,A \rangle A^i e_0 \wedge e_i +u^{-2} \langle U,A \rangle e_0 \wedge e_N -
\\ & 
\frac12 u^{-4} \langle A,A \rangle U^i e_0 \wedge  e_i 
+\frac12 u^{-2} k U^i e_0 \wedge e_i+u^{-2}U^iA^j e_i \wedge e_j-U^i e_i \wedge e_N.
\end{align*}
To put the basis in lexicographic order, we have to anti-symmetrise the term $U^iA^j e_i \wedge e_j$ to get $i<j$. By sorting and grouping the same terms, we obtain the expression of the statement.
 \end{proof}


\subsection{Conserved quantities} \label{subsec:conservedquant}   
We construct here conserved quantities of conformal circles on the sphere.

\begin{prop} \label{prop_conservedqtractor}
Consider the basis $b_2$ \eqref{baze} of $\wedge^2\R^{n+2}$. Pairing the $2$-tractor $\mathbb{Y}$ from  \eqref{2tractor} and the parallel tractors corresponding to basis elements $e_{ab}$ with respect to the tractor metric $G$ \eqref{met}, determines the conserved quantities $C^{ab}$ of conformal circles with constant $k$ as follows;
\begin{align} \label{CiN}
C^{iN} &= -(2u^{-6} \langle U,A \rangle^2   -
\frac12 u^{-4} \langle A,A \rangle   +\frac12 u^{-2} k )U^i + u^{-4} \langle U,A \rangle A^i  \\
\label{C0N}
C^{0N}&=-u^{-2}\langle U,A \rangle-\langle X, C^N \rangle\\ \label{Cij}
C^{ij} &= 
 u^{-2} \epsilon^{ij}(U, A) - \epsilon^{ij}(X, C^{N}) \\ 
\begin{split} \label{C0i}
C^{0i} &= U^i-u^{-2} \langle U,X \rangle A^i+u^{-2} \langle X,A \rangle U^i- u^{-2} \langle U,A \rangle X^i \\ &+\frac12 \|X\|^{2} C^{iN}-\langle X,C^N\rangle X^i
\end{split}
\end{align}
Here, $C^N$ is the vector of quantities
$$C^N=(C^{1N}, \dots, C^{nN})=(C^{iN}), \ i=1,\dots, n. $$
\end{prop}
\begin{proof}
Let us first describe basis parallel tractors by means of $\wedge^2\R^{n+2}$ in exponential coordinates. Thus we compute
 \begin{gather} \label{expw}
 \exp(-\rho_2(J))(w) = w-\rho_2(J)(w)+\frac12  \rho_2^2(J)(w), \ \ \ \ J=\left( \begin{smallmatrix} 0&0&0 \\ X &0&0 \\ 0&-X^t&0 \end{smallmatrix} \right) 
 \end{gather}
 for the representation \eqref{repr} for some neighborhood $X$ of $0$ in $\fc$, and $w=(w^{ab})$, $0\leq a <b \leq N$ are coordinates in the representation space $\wedge^2 \R^{n+2}$.
The series stops as $\rho_2^3(J)$ vanishes. 

The standard representation $\frak{so}(n+1,1) \to \frak{gl}(\R^{n+2})$ takes, as a linear map prescribed on the standard basis $e_a$, $a=0, \dots, N$, the form 
\begin{gather} \label{action2}
e_0 \mapsto \sum_{i=1}^n X^ie_i, \ \ \ \ \  e_i \mapsto -X^ie_{N}, \ \ \ \ \   e_N \mapsto 0, \ \ \ \ \   X=(X^i),\ i=1, \dots, n.
\end{gather}
Using \eqref{action2} and tensoriality, we compute the action $\rho_2$ explicitly. Since the standard action differs for $e_0$, $e_i$ and $e_N$, we have to consider four different types of basis elements of $b_2$. We get the action $\rho_2$ written, as a linear transformation prescribed on the basis $b_2$, as follows
\begin{align*} 
e_{0i} &\mapsto \sum_{j<i}X^j e_{ji}-\sum_{j>i}X^j e_{ij}-X^i e_{0N} \\
e_{0N} & \mapsto \sum_{j=1}^n X^j e_{jN} \\
e_{iN} &\mapsto 0\\
e_{ij} & \mapsto X^i e_{jN}-X^j e_{iN}
\end{align*}
Then we get the action $\rho_2^2$ as follows
\begin{align*}
e_{0i} &\mapsto \| X \|^2 e_{iN}-2X^i  \sum_{j=1}^n X^j e_{jN}  
\\
e_{0N} & \mapsto 0 \\
e_{iN} &\mapsto 0\\
e_{ij} & \mapsto 0 
\end{align*}
Altogether, we get according to \eqref{expw}
\begin{align*}
\exp(-\rho_2(J))(w)&= \sum_{a<b }w^{ab}e_{ab}-\Big(  \sum_{i=1}^n w^{0i}\big( \sum_{j<i}X^j e_{ji}-\sum_{j>i}X^j e_{ij}-X^i e_{0N}  \big) \\ &+w^{0N} \sum_{j=1}^n X^j e_{jN} + \sum_{i<j}^n w^{ij}(X^i e_{jN}-X^j e_{iN})   \Big)
\\ &+\frac12 \sum_{i=1}^n w^{0i}\big( \| X \|^2 e_{iN}-2X^i  \sum_{j=1}^n X^j e_{jN} \big),
\end{align*}
and after sorting the term according to basis elements, we get the right-hand side of the above equation as
\begin{align*}
& \quad \sum_{i=1}^n w^{0i} e_{0i} + \Big( w^{0N}+ \sum_{i=1}^n w^{0i}X^i  \Big)e_{0N}  \\ &+\sum_{i=1}^n \Big( w^{iN}-X^{i}w^{0N} +\sum_{j>i}X^{j}w^{ij}-\sum_{j<i}X^{j}w^{ji} +\frac{1}{2}\lVert X\rVert^{2}w^{0i} -X^{i}\sum_{j=1}^{n}X^{j}w^{0j}  \Big)e_{iN}  \\ &+ \sum_{i <j}^{n} \; \Big( w^{ij}-X^{j}w^{0i}+X^{i}w^{0j} \Big) e_{ij}.
\end{align*}
for $a,b=0, \dots, N$ and $i,j=1, \dots,n$.
For each basis element $e_{ab}$, i.e. for each coordinate $w^{ab}$, $0 \leq a <b \leq N$, we always choose one of the coordinates $w^{ab}$ to be $1$ and the rest to be $0$ in $\exp(-\rho_2(J))(w)$.
The pairing of the resulting expressions and the $2$-tractor $\mathbb{Y}$ from \eqref{2tractor} with respect to the tractor metric $G$ then means \emph{swapping $0$ and $N$ in one of the two expressions and standard pairing.} Doing these tedious but direct pairings, we get exactly the formulae from the statement.
\end{proof}

\section{Comparison of conserved quantities}
 We show here that the quantities of conformal circles discussed in previous sections are related. We formulate here the statement on comparison of both types of quantities; those coming from Noether’s principle applied to Conformal Mercator equation in Proposition \ref{prop_FTFRFDFS} and those coming from tractor calculations in Proposition \ref{prop_conservedqtractor}.

\subsection{Comparison of two types of quantities} \label{subsec:compare}
Consider conformal circles such that $\alpha_1'=0$. Under this assumption, they solve the Conformal Mercator equation \eqref{confMeqflat}. We know from Section \ref{sec:MercSol} that this assumption in not restrictive for traces of the curves.

\begin{prop} \label{prop_corrquantcircle}
The conserved quantities along conformal circles of Proposition \ref{prop_conservedqtractor} that are constructed from tractors are variational in the sense that they can be identified with the conserved quantities derived from Noether's principle of Proposition \ref{prop_FTFRFDFS}, with $C$ given by \eqref{Ck}, by the following relations:
\begin{align} 
C^{iN} &= -E_T^ i \label{CiNET}\\ 
C^{0N}&= E_D \label{C0NED} \\ 
C^{ij} &=  E_R^{ij}\\ 
C^{0i} &= -\frac{1}{2} E_S^i. \label{C0iES}
\end{align}
\end{prop}
\begin{proof}
The vector  $C^N = (C^{iN})$ of Proposition \ref{prop_conservedqtractor}
 is precisely the vector $C$ \eqref{Ck}. Thus we have the identity \eqref{CiNET} for conformal circles with constant $\a_1=k.$  The rest of the identities, \eqref{C0NED} - \eqref{C0iES}, follow immediately.
\end{proof}

\subsection{Examples} 
We give here explicit examples of conserved quantities; we follow  the tractor methods to explain the details on examples.

\begin{exam}{\it{Quantities in dimension $3$.}} \label{ex1}
We have $10$ basis elements of the lexicographic basis of $\wedge^2(\R^5)$
$$e_1=e_{01}, e_2=e_{02}, e_3=e_{03},e_4=e_{04},  e_5=e_{12}, e_6=e_{13}, e_7=e_{14}, e_8=e_{23}, e_9=e_{24}, e_{10}=e_{34}$$ and we denote corresponding coordinates $w=(w^1, \dots, w^{10})^t=(w^a)^t$, $a=1,\dots,10$. We find $\rho_2(J)$ as a linear map prescribed on the basis $(e_i), i=1\dots 10$ as the following  transformation matrix 
$$
\rho_2(J)= \left[ \begin {smallmatrix} 0&0&0&0&0&0&0&0&0&0
\\ 0&0&0&0&0&0&0&0&0&0\\ 0&0&0&0&0
&0&0&0&0&0\\ -{ X_1}&-{ X_2}&-{ X_3}&0&0&0&0&0
&0&0\\ -{ X_2}&{ X_1}&0&0&0&0&0&0&0&0
\\ -{ X_3}&0&{ X_1}&0&0&0&0&0&0&0
\\ 0&0&0&{ X_1}&-{ X_2}&-{ X_3}&0&0&0&0
\\ 0&-{ X_3}&{ X_2}&0&0&0&0&0&0&0
\\ 0&0&0&{ X_2}&{ X_1}&0&0&-{ X_3}&0&0
\\ 0&0&0&{ X_3}&0&{ X_1}&0&{ X_2}&0&0
\end {smallmatrix} \right].
$$
Note that \emph{same types} of basis elements are next to each other, with the exception $e_7$ and $e_8$, one can swap $7$th and $8$th row and column to see the representation \emph{more logically}.
Then the linear transformation $\exp(-\rho_2(J))$ takes in the above basis form
$$
\exp(-\rho_2(J))= \left[ \begin{smallmatrix} 1&0&0&0&0&0&0&0&0&0
\\ 0&1&0&0&0&0&0&0&0&0\\ 0&0&1&0&0
&0&0&0&0&0\\ { X_1}&{ X_2}&{ X_3}&1&0&0&0&0&0&0
\\ { X_2}&-{ X_1}&0&0&1&0&0&0&0&0
\\ { X_3}&0&-{ X_1}&0&0&1&0&0&0&0
\\ \frac12(- X_1^{2}+ X_2^{2}+
 X_3^{2})&-{ X_1}{ X_2}&-{ X_1}{ X_3}&-{ X_1}&{ X_2
}&{ X_3}&1&0&0&0\\ 0&{ X_3}&-{ X_2}&0&0&0&0&1&0
&0\\ -{ X_1}{ X_2}&\frac12( X_1^{2}-
 X_2^{2}+ X_3^{2})&-{ X_2}{ X_3}&-{ X_2}&-{ 
X_1}&0&0&{ X_3}&1&0\\ -{ X_1}{ X_3}&-{ X_2}
{ X_3}&\frac12( X_1^{2}+X_2^{2}-X_3^{2})
&-{ X_3}&0&-{ X_1}&0&-{ X_2}&0&1 \end{smallmatrix} \right] 
$$
Then one describes elements of $\wedge^2\R^5$ in exponential coordinates easily via matrix multiplication as $\exp(-\rho_2(J))\cdot w$.

Now, pairing these with $\mathbb{Y}$ with respect to the tractor metric gives conserved quantities. We obtain a generic family of quantities so that for each $k=1, \dots, 10$, we choose $w^k=1$ and the remaining coordinates zero.
Altogether, we get $10$ quantities as follows, where $X=(x_1,x_2,x_3)$, $U=(\dot x_1, \dot x_2, \dot x_3)$, $A=(\ddot x_1, \ddot x_2, \ddot x_3)$ and $A'=(\dddot x_1, \dddot x_2, \dddot x_3)$.
We have $3$ quantities of type \eqref{CiN} corresponding to $e_7, e_9, e_{10}$. We write them as a vector $C^{i4}:=(C^{14},C^{24},C^{34})$
\begin{align*}
C^{i4}={\frac { \langle U, A \rangle ^{2}U}{ \langle U, U \rangle ^{3}}}-{\frac { ( 6u^{-2}\langle U,A \rangle^2- \langle A,A \rangle +u^2k) U}{ 2\langle U, U \rangle ^{2}}}+{\frac { \langle U, A \rangle A}{ \langle U, U \rangle ^{2}}}
\end{align*}
We have $1$ quantity of type \eqref{C0N} corresponding to $e_4$ of the form
\begin{align*}
C^{04}= -{\frac {\langle U, A \rangle }{\langle U, U \rangle}}-{\frac { \langle U, A \rangle ^{2} \langle X, U \rangle }{
 \langle U, U \rangle ^{3}}}+{
\frac { ( 6u^{-2}\langle U,A \rangle^2- \langle A,A \rangle +u^2k )  \langle X, U \rangle }{ 2\langle U, U \rangle ^{2}}}-{\frac {
 \langle U, A \rangle  \langle X, A \rangle }{ \langle U, U \rangle ^{2}}}
\end{align*}
We have $3$ quantities of type \eqref{Cij} corresponding to $e_5$, $e_6$ and $e_8$. We write them as a vector $C^{ij}:=(C^{12},-C^{13},C^{23})$ as follows, where $V_1 \wedge  V_2$ denotes vector product of vectors $V_1, V_2$
 \begin{align*}
C^{ij}={\frac {1}{\langle U, U \rangle}U \wedge A}-{\frac { \langle U, A \rangle ^{2}  }{ \langle U, U \rangle ^{3}} X \wedge U}+{\frac {  6u^{-2}\langle U,A \rangle^2- \langle A,A \rangle +u^2k  }{
 2\langle U, U \rangle^{2}}} X \wedge U-{
\frac { \langle U, A \rangle   }{
 \langle U, U \rangle^{2}}} X \wedge A.
 \end{align*}
We have $3$ quantities of type \eqref{C0i} corresponding to $e_1,e_2,e_3$ written as a vector $C^{0i}:=(C^{01}, C^{02}, C^{03})$ 
\begin{align*}
C^{0i} &=
U+\frac{\langle X,A \rangle}{\langle U,U \rangle}U-\frac{\langle X,U \rangle}{\langle U,U \rangle}A- \frac{\langle U,A \rangle}{\langle U,U \rangle}X
 \\ &-\frac12 \langle X, X \rangle  ( -{\frac {  \langle U, A \rangle  ^{2}U}{ \langle U, U \rangle ^{3}}}+{\frac { ( 
6u^{-2}\langle U,A \rangle^2-\langle A,A \rangle +u^2k ) U}{ 2\langle U, U \rangle  ^{2}}}-{\frac { \langle U, A \rangle  A}{
 \langle U, U \rangle^{2}}}
 ) \\ &+( -{\frac { \langle U, A \rangle  ^{2} \langle X, U \rangle  }{ \langle U, U \rangle^{3}}}+{\frac { ( 6u^{-2}\langle U,A \rangle^2- \langle A,A \rangle +u^2k )  \langle X, U \rangle  }{ 2\langle U, U \rangle  ^{2}}}-{\frac { \langle U, A \rangle   \langle X, A \rangle  }{ \langle U, U \rangle  ^{2}}} ) X
\end{align*}

\end{exam}

\begin{exam}{\it{Quantities for projectively parametrized conformal circles in dimension $3$.}}
Let us consider the open family of conformal circles   parametrized by three parametres $a, c_1, c_2$ as
\begin{align*}
x_1(t)&= { p_1}+\frac {\sin(c_2)(a\cos(c_1)t^2+\sin(c_1)t)}{a^2t^2+1}\\
  x_2(t)&= { p_2}+{\frac {\cos(c_2)(a\cos(c_1)t^2+\sin(c_1)t)}{a^2t^2+1}} \\
  x_3(t)&= { p_3}+{\frac {-a\sin ( { c_1} ) {t}^{2}+\cos ( {
 c_1} ) t}{a^2t^2+1 }}
 \end{align*}
 For $a=0$, we get a line, otherwise we get a circle of radius $R=\frac{1}{2|a|}$. The circle lies in the plane given by the point $P=(p_1,p_2,p_3)$ and directions $(\sin(c_2),\cos(c_2),0)$ and $(0,0,1)$. Direct computation of derivatives gives
 \begin{gather*}
u^2=\langle U,U \rangle =\frac{1}{(a^2t^2+1)^2}, \ \langle U,A \rangle =-\frac{2a^2t}{(a^2t^2+1)^3}, \\ 
\langle A,A \rangle =\frac{4a^2}{(a^2t^2+1)^3}, \  
\langle A,A' \rangle =-\frac{12a^4t}{(a^2t^2+1)^4}, \ \langle A',A' \rangle =\frac{36a^4}{(a^2t^2+1)^4}
 \end{gather*}
 
Direct computation with these expressions
confirms that $\Delta_4=0$ and we also get $\alpha_1=0$, i.e. we get conformal circles that indeed carry projective parametrisation. We analogously compute $\langle X,X \rangle$, $\langle X,U \rangle $, $\langle X,A \rangle $, $\langle X,A' \rangle $ (that are nastier, so we do not give them explicitly; we use Maple). Substituting all these products into the quantities from the previous example gives the quantities.

The quantities corresponding to $e_i$, $i=1, \dots, 10$ (in this order) respectively take the form
\begin{align*} 
{C}_{{1}}&= (  ( (- p_1^{2}+p_2^{2}+
 p_3^{2})a^2+1 ) \sin ( { c_2} ) -2a^2{ p_1}{
 p_2}\cos ( { c_2} )  ) \sin ( { c_1}
 ) -2a{ p_3} ( \cos ( { c_1} ) a{ p_1}
+\sin ( { c_2} )  ) \\
{C}_{{2}}&= (  (  (p_1^{2}- p_2^{2}+
 p_3^{2})a^2+1 ) \sin ( { c_1} ) -2a{ p_3}
 ) \cos ( { c_2} ) -2a^2{ p_2} ( \sin
 ( { c_2} ) \sin ( { c_1} ) { p_1}+\cos
 ( { c_1} ) { p_3} ) \\
{C}_{{3}}&= (  (p_1^{2}+ p_2^{2}- p_3^
{2})a^2+1 ) \cos ( { c_1} ) -2a(\sin(c_2)p_1+\cos(c_2)p_2)(ap_3\sin(c_1)-1) \\
{C}_{{4}}&=-2a^2(\sin ( { c_2} ) \sin ( { c_1
} ) { p_1}+\sin ( { c_1} ) \cos ( { 
c_2} ) { p_2}+\cos ( { c_1} ) { p_3} )\\
{C}_{{5}}&=2a^2 ( \sin ( { c_2} ) { p_2}-
\cos ( { c_2} ) { p_1} ) \sin ( { c_1}
 ) \\
{C}_{{6}}&=2a(\sin ( { c_2} ) \sin ( { c_1}
 ) a{ p_3}-\cos ( { c_1} ) a{ p_1}-\sin
 ( { c_2} ) )\\
{C}_{{7}}&=2a^2\sin ( { c_1} ) \sin ( { c_2}
 ) \\
{C}_{{8}}&=2a(\sin ( { c_1} ) \cos ( { c_2}
 ) a{ p_3}-\cos ( { c_1} )a { p_2}-\cos
 ( { c_2} )) \\
{C}_{{9}}&=2a^2\sin ( { c_1} ) \cos ( { c_2}
 ) \\
{C}_{{10}}&=2a^2\cos ( { c_1} ) 
\end{align*}

One can check by direct calculations that $$C_7^2+C_9^2+C_{10}^2=4a^2=\frac{1}{2R^2}.$$ In fact, denoting 
$T=(C_7,C_9,C_{10}),$ the initial velocity at the origin takes form $U(0)=R^2T$. 
Moreover, $$-C_4=\langle P, T \rangle =p_1C_7+p_2C_9+p_3C_{10}.$$
In fact, the quantities $C_7,C_9,C_{10}$ correspond to quantities derived from the generators of translation and $C_4$ corresponds to that from the generator of dilatation according to Propositions \ref{prop_corrquantcircle} and \ref{prop_FTFRFDFS}.
We also compute 
\begin{align*}
C_5&=p_2C_7-p_1C_9\\
C_6&=p_3C_7-p_1C_{10}-2a \sin(c_2)\\
C_8&=p_3C_9-p_2C_{10}-2a\cos(c_2),
\end{align*}
which come from the generators of rotation.
Then
$$(C_6-p_3C_7+p_1C_{10})^2+(C_8-p_3C_9+p_2C_{10})^2=4a^2.$$
According to all these observations, we can also write the quantities constructed from the generators of special conformal transformation, $C_1, C_2, C_3,$ as
$$
(C_1, C_2, C_3)  = \left(1+a^2|P|^2\right)U(0) +C_4 P -2a (P \w n),
$$
where $n=(\cos(c_2), -\sin(c_2), 0)$ is a unit normal vector to the vector $U(0)$.

\end{exam}


\section*{Acknowledgements}
This article is based upon work from COST Action CaLISTA CA21109 supported by COST. The authors would like to thank Jan Gregorovi\v c, Josef \v Silhan and Vojt\v ech \v Z\' adn\' ik for helpful discussions based on \cite{krivule}.  
P.P. is grateful to the Department of Mathematics, Faculty of Science, University of South Bohemia for kind hospitality and support during her research visit.

AI disclosure: Copilot was used by L.Z. to check and translate long expressions from CAS Maple to LaTeX.

\end{document}